\documentclass[10pt,reqno]{amsart}
\usepackage[cp1251]{inputenc}
\usepackage[russian]{babel}
\usepackage{amsmath}
\usepackage{amssymb}
\usepackage{amsfonts}
\usepackage{graphicx}
\usepackage{eucal}
\usepackage{hyperref,amsthm}
\RequirePackage{bbold}

\newtheorem{theorem}{Теорема}

\theoremstyle{definition}

\numberwithin{equation}{section}
\begin{document}

\title[Об оценке структурного параметра Марковского Q-процесса]
    {Об оценке структурного параметра Марковского Q-процесса}

\author{\large{Аъзам А.~Имомов, \, Зухриддин А.~Назаров}}
\address {Аъзам Абдурахимович Имомов, Зухриддин Назаров
\newline\hphantom{iii} Каршинский государственный университет,
\newline\hphantom{iii} Институт Математики имени В.И.Романовского, Узбекистан.}
\email{{{imomov{\_}\,azam@mail.ru}, {zuhrov13@gmail.com} }}


\subjclass[2010] {Primary 60J80; Secondary 60J85}



\begin{abstract}
    In the paper we consider a stochastic model which called Markov Q-processes that forms a continuous-time Markov population system.
    Markov Q-processes are defined as stochastic Markov branching processes with trajectories continuing in the remote future.
    Estimation of the structural parameter of the Markov Q-process is the main goal of this paper. To estimate this parameter,
    an unbiased estimator of the Lotka-Nagaev type is proposed. An asymptotic expansion of the variance of this estimator is found.

\vspace{3mm}

\noindent {\textsc{Ключевые слова и фразы}}. \textit{Марковские ветвящиеся системы; Марковские Q-процессы;
    переходные вероятности; производящая функция; структурный параметр; несмещенная оценка}.
\end{abstract}

\maketitle

\section{Введение и результаты}    \label{MySec:1}

    Модели ветвящихся случайных систем являются существенной частью общей теории случайных процессов. Растущий интерес к этим
    моделям обусловлен многими факторами. Первой из этих причин, которая стала основным толчком к созданию теории ветвящихся
    случайных моделей, является возможность оценить с их помощью вероятности выживания популяции однотипных индивидуумов;
    см. \cite{Sev71}, \cite{ANey} и \cite{Harris}. Среди всевозможных случайных траекторий всех моделей ветвящихся систем есть
    такие, которые длятся бесконечно долго. В случае модели Гальтона-Ватсона, класс траекторий не вырождающихся в далеком будущем,
    образует так называемые Q-процессы; см. \cite{ANey} и \cite{Imomov2014}. В случае марковских ветвящихся случайных систем с
    непрерывным временем аналогичная модель называемая \emph{марковский Q-процесс}, впервые определена в работе \cite{Imomov2012}.

    В работе \cite{Imomov2005} доказан дифференциальный аналог основной леммы теории Марковских ветвящихся случайных систем
    непрерывного времени. С помощью этой леммы в работе \cite{Imomov2012} исследованы структурные свойства марковских Q-процессов.
    В работе \cite{Imomov2014} установлена глубокая связь между Q-процессами и ветвящимися системами Гальтона-Ватсона с иммиграцией.
    Там же доказана предельная теорема для совместных распределений состояний и общих состояний в Q-процессе.

    В данной работе нас интересует задача оценки структурного параметра марковского Q-процесса. Для оценки этого параметра
    предлагается несмещенная оценка типа Лотка-Нагаева и исследуются асимптотические поведения ее дисперсии.

    Обозначим через ${\mathbb{N}}$ множество натуральных чисел и пусть ${\mathbb{N}}_0=\{0\}\cup{\mathbb{N}}$. Случайной величиной
    $Z(t)$ обозначим размер популяции в момент $t\geq 0$ в однородно-непрерывной во времени Марковской ветвящейся случайной
    системе (МВС) с интенсивностью ветвления $\left\{ {a_k,\;k \in{\mathbb{N}}_0}\right\}$. Каждая частица в этой системе
    имеет экспоненциально распределенный случайный период жизни со средним значением $\sum_{k\neq{1}}{a_k}$ и она в конце
    своей жизни производит $k\in{\mathbb{N}}_0 \backslash \{1\}$ потомков с вероятностью ${-a_{k}}\left/ {a_{1}}\right.$.
    Соответствующая интенсивности ветвления \emph{q-матрица}, ${\textbf{Q}}=\left\{q_{ij}\right\}$, с компонентами
\begin{eqnarray}           \label{eq1}
    q_{ij}
    = \left\{\begin{array}{l} ia_{0} \;,   \qquad\qquad \hfill\text{если} \quad j = i-1,   \\
\\
    ia_{j-i+1} \;,   \qquad\qquad \hfill\text{если} \quad j \geq i \geq 0, \\
\\
    0 \;,  \qquad\qquad\hfill \text{в остальных случаях}, \\
    \end{array} \right.
\end{eqnarray}
     где
\begin{equation*}
    a_k\ge 0 \quad \text{для} \quad k\in{\mathbb{N}}_0 \backslash \{1\} \qquad \text{и} \qquad
    0< a_0<-a_1=\sum_{k\in{\mathbb{N}}_0 \backslash \{1\}}{a_k} < \infty,
\end{equation*}
     полностью характеризует дальнейшую эволюцию МВС; см. \cite{LiDEDS21}. Определенная выше система образует разложимую и
     однородно-непрерывную во времени цепь Маркова с пространством состояний, состоящим из двух классов:
     $\mathcal{S}_{0}=\left\{0\right\}\cup{\mathcal{S}}$, здесь $\mathcal{S}\subset{\mathbb{N}}$.  При этом состояние
     $\left\{0\right\}$ является поглощающим, а $\mathcal{S}$ -- класс возможных существенных сообщающихся состояний.

     Введем в рассмотрение переходные вероятности
\begin{equation*}
    P_{ij}(t):={\textsf{P}}\Bigl\{{Z(t+\tau)=j\bigm|{Z(\tau)=i}}\Bigr\} \quad \text{для любого} \quad \tau \ge 0.
\end{equation*}
    Эти вероятности равны $i$-кратной свертке распределения $P_{1j}(t)$, т.е.
\begin{equation*}
    P_{ij}(t)=\sum\limits_{j_1+j_2+\ldots+j_i=j}{P_{1j_1}(t)\cdot P_{1j_2}(t)\cdot \, \ldots \, \cdot P_{1j_i}(t)}.
\end{equation*}
    В свою очередь, с помощью q-матрицы \eqref{eq1} можно вычислить, что вероятности
    $P_{1j}(t)$ допускают следующее локальное представление:
\begin{equation}\label{eq2}
    P_{1j}(\varepsilon)=\delta_{1j}+a_j\varepsilon+o(\varepsilon) \quad \text{при} \quad \varepsilon\downarrow 0,
\end{equation}
    где $\delta_{1j}$ -- знак Кронекера; см.~\cite{Sev71}.

    Пусть $\textsf{P}_{i}\bigl\{{*}\bigr\}:= \textsf{P}\left\{{{\ast} \bigm| {Z(0)=i}}\right\}$ и рассмотрим условные вероятности
\begin{equation*}
    \textsf{P}_i^{\mathcal{H}(t+u)}\{*\}: = \textsf{P}_i
    \left\{{*\bigm| {t+u<{\mathcal{H}}< \infty}} \right\} \quad \text{для некоторого} \quad  u\geq 0,
\end{equation*}
    где величина $\mathcal{H} := \inf\bigl\{t: Z(t)=0\bigr\}$ обозначает момент вырождения МВС.
    В работе {\cite{Imomov2012}} доказано, что
\begin{equation}               \label{eq3}
    \mathcal{Q}_{ij}(t): = \mathop {\lim }\limits_{u \to \infty}
    \textsf{P}_i^{{\mathcal{H}}(t+u)}\bigl\{Z(t)=j\bigr\}={{jq^{j-i}}\over{i\beta^t}}P_{ij}(t),
\end{equation}
    где
\begin{equation*}
    \beta:=\exp\left\{{\sum_{j\in{\mathbb{N}}}{ja_j q^{j-1}}}\right\},
\end{equation*}
    а число $q$ -- вероятность исчезновения одной частицы претерпевающей превращения по закону интенсивности
    $\left\{ {a_k,\;k \in{\mathbb{N}}_0}\right\}$ такая, что 
\begin{equation*}
    q=\inf\bigl\{x\in[0,1]: \; f(x)=0\bigr\}.
\end{equation*}
    Нетрудно проверить, что ${\sum_{j\in{\mathbb{N}}}{\mathcal{Q}_{ij}(t)}}=1$ при каждом $i\in{\mathbb{N}}$.
    Вероятностная мера $\bigl\{\mathcal{Q}_{ij}(t)\bigr\}$ определяет новый процесс развития популяционной
    системы -- непрерывно-однородную неразложимую марковскую цепь с пространством всевозможных состояний
    ${\mathcal{E}} \subset \mathbb{N}$, называемый \textit{Марковский Q-процесс} (МQП)
    $\left\{{W(t)}, t\geq{0}\right\}$. По определению
\begin{equation*}
    {\mathcal{Q}}_{ij}(t) = \textsf{P}\left\{W(t)=j \bigm| {W(0)=i}\right\}
    = \textsf{P}_i \left\{{Z(t)=j \bigm| {{\mathcal{H}}=\infty}}\right\},
\end{equation*}
    так что МQП можно интерпретировать как ``долгоживущая'' МВС; см.~{\cite{Imomov2022}}. Отсюда, с помощью \eqref{eq2},
    находим следующее локальное представление вероятностей $\mathcal{Q}_{1j}(\varepsilon)$:
\begin{equation*}
    \mathcal{Q}_{1j}(\varepsilon)=\delta_{1j}+p_j\varepsilon+o(\varepsilon ) \quad \text{при} \quad \varepsilon\downarrow 0,
\end{equation*}
    с плотностями вероятностей перехода
\begin{equation*}
    p_0=0, \quad p_1=a_1-\ln\beta<0 \quad \text{и} \quad p_j=jq^{j-1}a_j\ge 0
    \quad \text{для} \quad  j\in{\mathcal{E}}\backslash\{1\}.
\end{equation*}
    Следовательно, производящая функция (ПФ)
\begin{equation*}
    g(x):=\sum_{j\in{\mathcal{E}}}{p_j x^j}=x\bigl[{f'(qx)-f'(q)}\bigr]
\end{equation*}
    полностью определяет МQП, здесь $f(x)$ есть инфинитезимальная ПФ, порождающая МВС $\left\{{Z(t)}, t\geq{0} \right\}$,
    то есть $f(x)=\sum_{j\in{\mathcal{S}}_0}{a_j x^j}$. В этом обозначении (см.~\cite{Imomov2012})
\begin{equation*}
    \beta=\exp\bigl\{{f'(q)}\bigr\} \quad \text{и} \quad f(q)=0.
\end{equation*}

    Известно, что регулирующим параметром для МВС вступает величина $f'(1-)$ и по свойству траекторий выделяется три типа процессов,
    характеризующихся его значением. МВС называется докритическим, если $f'(1-)<1$, критическим, если $f'(1-)=1$, и надкритическим,
    если $f'(1-)>1$. Соответственно, эволюция МQП управляется (регулируется) по существу структурным параметром $\beta$ и
    известно, что $\beta=1$ при $f'(1-)=0$ и $\beta<1$ если $f'(1-)\ne 0$; см. \cite{Imomov2012} и литературу в ней.
    В соответствие с случаем дискретного времени \cite[с.~59,~Теорема~2]{ANey}, ${\mathcal{E}}$ является положительно возвратным
    если $\beta <1$ и, ${\mathcal{E}}$ является невозвратным, если $\beta =1$.
    Таким образом, различается два типа МQП в зависимости от значений параметра $\beta$. Сказанное наглядно
    утверждается в следующей предельной теореме, по предельным свойствам локальной вероятности $\mathcal{Q}_{11}(t)$.

\medskip

\textbf{Теорема~A}~\cite{Imomov2012}. \emph{Пусть имеется МQП $\left\{ {W(t)}, t\geq{0} \right\}$,
                    порожденный инфинитезимальной ПФ $g(x)$ с первым конечным моментом $b:=g'(1-)<\infty$.}
\begin{itemize}
    \item    \emph{Если $\beta=1$, то
\begin{equation}               \label{eq4}
    t^2 \mathcal{Q}_{11}(t)\longrightarrow{2\over{ba_0}} \quad \text{при} \quad t\to\infty.
\end{equation}}
    \item    \emph{Если $\beta<1$, то
\begin{equation}              \label{eq5}
    \mathcal{Q}_{11}(t)\longrightarrow{{\left|{\ln\beta}\right|{\mathcal{A}}}\over{a_0}} \quad \text{при} \quad t\to\infty,
\end{equation}
    где постоянная
\begin{equation*}
    \mathcal{A}=q\exp\left\{{\int_0^q{\left[{{1\over{s-q}}-{{f'(q)}\over{f(s)}}}\right]ds}}\right\}<\infty.
\end{equation*}}
\end{itemize}

\medskip

    Отметим, что при $q =1$, положительная постоянная $\mathcal{A}$ представляет собой
    константу Колмогорова-Севастьянова из теории докритических МВС ($f'(1-)<0$); см. \cite{Sev71}.

Введем теперь в рассмотрение ПФ распределения состояний МQП
\begin{equation*}
    G_i(t;x):={\textsf{E}}_i x^{W(t)}={\textsf{E}}\left[ {{x^{W(t)}}\Bigm|W(0)=i}\right]
    =\sum_{j\in{\mathcal{E}}}{\mathcal{Q}_{ij}(t)x^j}.
\end{equation*}
    Как было доказано в \cite{Imomov2012},
\begin{equation}              \label{eq6}
    G_i (t;x)=x\left[{{{\Phi(t;qx)}\over q}}\right]^{i-1}
    \exp\left\{{\int_0^t{b\left({{{\Phi(\tau;qx)}\over q}}\right)d\tau}}\right\},
\end{equation}
    где $\Phi(t;x)=\sum\nolimits_{j\in{\mathcal{S}}_0}{P_{1j}(t)x^j}$ и
\begin{equation*}
    b(x)={{g(x)}\over x}=f'(qx)-f'(q).
\end{equation*}
    Отсюда, путем дифференцирования в точке $x =1$ из получаем
\begin{eqnarray}           \label{eq7}
    {\textsf{E}}_i W(t)=\bigl({i-1}\bigr)\beta^t+{\textsf{E}}_1 W(t)
    = \left\{\begin{array}{l} \bigl({i-1}\bigr)\beta^t+bt+1   \hfill, \quad \text{если} \quad \beta = 1,  \\
\\
    \bigl({i-1}\bigr)\beta^t+1+\gamma\bigl(1-\beta^t\bigr) \; \hfill, \quad \text{если} \quad \beta < 1,
    \end{array} \right.
\end{eqnarray}
    и
\begin{eqnarray}           \label{eq8}
    {\textsf{Var}}_i W(t)
    = \left\{\begin{array}{l} bti   \hfill, \quad \text{если} \quad \beta = 1,  \\
\\
    \Bigl[{\gamma+\left({i-1}\right)\left({1+\gamma}\right)\beta^t}\Bigr]
    \bigl({1-\beta^t}\bigr) \; \hfill, \quad \text{если} \quad \beta < 1,
    \end{array} \right.
\end{eqnarray}
    где
\begin{equation*}
    {\textsf{Var}}_i W(t)={\textsf{Var}}\left[{W(t)\bigm|{W(0)=i}}\right]
    \quad \text{и} \quad  \gamma:={b\big/ \bigl|{\ln\beta}\bigr|}.
\end{equation*}

    Ввиду всего вышеизложенного, поставим задачу оценить параметр $\beta$
    по наблюденным значениям $W(t)$. Из формулы \eqref{eq6} следует, что
\begin{equation*}
    {\textsf{E}}\left[{{x^{W(t+1)}}\bigm|W(t)}\right]
    =\left[{{{\Phi(1;qx)}\over q}}\right]^{W(t)-1} G(1;x),
\end{equation*}
    где $G(t;x):=G_1(t;x)$. Отсюда получаем
\begin{equation*}
    {\textsf{E}}\left[{W(t+1)\bigm|{W(t)}}\right]=\bigl[{W(t)-1}\bigr] \beta+{\textsf{E}}_1 W(1).
\end{equation*}
    Последнее равенство позволит нам выписать уравнение
\begin{equation*}
    W(t+1)=\bigl[{W(t)-1}\bigr]\beta+W(1)+\epsilon(t),
\end{equation*}
    с погрешностью $\epsilon(t)$, имеющей нулевое среднее: ${\textsf{E}}\epsilon(t)=0$.
    Учитывая это уравнение, предлагаем следующую оценку для $\beta$, при известном ${\textsf{E}}_1 W(1)$:
\begin{equation*}
    \widehat\beta(t)={{W(t+1)-{\textsf{E}}_1 W(1)}\over{W(t)-1}} \raise0.9pt\hbox{,} \quad t>1.
\end{equation*}

    Оценка $\widehat\beta(t)$ является несмещенной для параметра $\beta$. Действительно,
    согласно формуле полной вероятности и однородности МQП, с учетом \eqref{eq7}
\begin{align*}
    {\textsf{E}}\widehat\beta(t)
    & = \sum_{j\in{\mathcal{E}}}{\textsf{P}\left\{{W(t)=j}\right\}{\textsf{E}}
    \left[{{{{W(t+1)-{\textsf{E}}_1 W(1)}\over{W(t)-1}}}\Bigm|{W(t)=j}}\right]} \\
\\
    & = \sum_{j\in{\mathcal{E}}} {\mathcal{Q}_{1j}(t){{\textsf{E}}_j W(1)-{\textsf{E}}_1 W(1)}\over{j-1}}
    =\beta\sum_{j\in{\mathcal{E}}}{\mathcal{Q}_{1j}(t)}=\beta.
\end{align*}

    Следующие теоремы характеризуют дальнейшие свойства оценки $\widehat\beta(t)$.

\begin{theorem}               \label{INTh:1}
    Пусть $b<\infty$. Если $\beta=1$, то
\begin{equation*}
    {t\over 2}\cdot{\textsf{Var}}\widehat\beta(t)=1+\mathcal{O}\left({{{\ln^2t}\over t}}\right)
    \quad \text{при}\quad t \to \infty.
\end{equation*}
\end{theorem}

\begin{theorem}               \label{INTh:2}
    Пусть $b<\infty$. Если $\beta<1$, то
\begin{equation*}
    {\textsf{Var}}\widehat\beta(t)=\mathcal{O}\left(1\right)  \quad \text{при}\quad t \to \infty.
\end{equation*}
\end{theorem}

    В следующем параграфе, наши рассуждения по части доказательства вышеизложенных
    теорем основываются на методе, использованном в работе А.В.Нагаева~\cite{Nagaev}.

\section{Доказательство Теорем 1 и 2}          \label{MySec:2}

    Согласно формуле полной вероятности, имеем
\begin{align*}
    {\textsf{Var}}\widehat\beta(t)
    & = {\textsf{E}}\left[\widehat\beta(t)-{\textsf{E}}\beta(t)\right]^{2}
    ={\textsf{E}}\left[\widehat\beta(t)-\beta\right]^{2} \\
\\
    & = \sum_{j\in{\mathcal{E}}}{\mathcal{Q}_{1j}(t){\textsf{E}}\left[ \left({{{W(t+1)-
    {\textsf{E}}_1 W(1)-\left(W(t)-1\right)\beta}\over{W(t)-1}}} \right)^{2}\Bigm|{W(t)=j}\right]}.
\end{align*}
    Далее, силу однородность МQП, получим
\begin{align}               \label{eq9}
    {\textsf{Var}}\widehat\beta(t)
    & = \sum_{j\in{\mathcal{E}\backslash \{1\}}}\mathcal{Q}_{1j}(t)\frac{1}{\bigl(j-1\bigr)^{2}}
    {\textsf{E}}_{j}\Bigl[W(1)-\textsf{E}_{1}W(1)-\bigl(W(0)-1\bigr)\beta\Bigr]^{2} \nonumber\\
\nonumber\\
    & = \sum_{j\in{\mathcal{E}\backslash \{1\}}}\mathcal{Q}_{1j}(t)\frac{1}{\bigl(j-1\bigr)^{2}}
    {\textsf{E}}_{j}\Bigl[W(1)-\textsf{E}_{1}W(1)\Bigr]^{2} \nonumber\\
\nonumber\\
    & = \sum_{j\in{\mathcal{E}\backslash \{1\}}}\mathcal{Q}_{1j}(t)\frac{1}{\bigl(j-1\bigr)^{2}}{\textsf{Var}}_{j}W(1)
        =\sum_{k\in{\mathcal{E}}}{{{\mathcal{Q}_{1k+1}(t)}\over{k^2}}{{\textsf{Var}}_{k+1} W(1)}}.
\end{align}

\begin{proof}[\textbf{Доказательство Теоремы~\ref{INTh:1}}]
    Из \eqref{eq9}, учитывая \eqref{eq8} для случая $\beta=1$, находим
\begin{align}               \label{eq10}
    {\textsf{Var}}\widehat\beta(t)
    & = \sum_{k\in{\mathcal{E}}}{{{\mathcal{Q}_{1k+1}(t)}\over{k^2}}{{\textsf{Var}}_{k+1} W(1)}}
    =\sum_{k\in{\mathcal{E}}}{{\mathcal{Q}_{1k+1}(t)}\over{k^2}} b(k+1)    \nonumber\\
\nonumber\\
    & = b\sum_{k\in{\mathcal{E}}}{{{\mathcal{Q}_{1k+1}(t)}\over{k}}}+
    b\sum_{k\in{\mathcal{E}}}{{{\mathcal{Q}_{1k+1}(t)}\over{k^2}}}=:\Sigma_1(t)+\Sigma_2(t).
\end{align}
    Согласно определению ПФ $G(t;x)$, с учетом того, что $\int_{0}^{1}x^{k-1}dx={1}\big/{k}$,
    первую сумму в \eqref{eq10} можно преобразовывать к следующему виду:
\begin{equation*}
    \Sigma_1(t)= b\sum_{k\in{\mathcal{E}}}{{{\mathcal{Q}_{1k+1}(t)}\over k}}
    = b\int_{0}^{1}\sum_{k\in{\mathcal{E}}}\mathcal{Q}_{1k+1}(t)x^{k-1}dx.
\end{equation*}
    Далее имеем
\begin{align}               \label{eq11}
    \Sigma_1(t)
    & = b\int_{0}^{1}\frac{1}{x^{2}}\left(\sum_{j\in{\mathcal{E}}}\mathcal{Q}_{1j}(t)x^{j}-
        \mathcal{Q}_{11}(t)x\right)dx   \nonumber\\
\nonumber\\
    & = b\int_0^1{{{G(t;x)-x\mathcal{Q}_{11}(t)}\over{x^2}}dx}=:b\int_0^1{{{\tau(t;x)}\over{x}} dx}.
\end{align}
    Очевидно, что $G(t;x)$ как вероятностная ПФ монотонно возрастает по $0\leq{x}<1$ и, следовательно,
    нетрудно проверить, что подынтегральная функция ${\tau(t;x)}\big/{x}$ в \eqref{eq11} также является
    монотонно возрастающей. Причем легко убедиться, что она имеет конечные значения в концах области интегрирования:
\begin{align*}
    \lim_{x\downarrow 0}{{\tau(t;x)}\over{x}} &=  \mathcal{Q}_{12}(t), \\
\nonumber\\
    \lim_{x\uparrow 1}{{\tau(t;x)}\over{x}} &=  1-\mathcal{Q}_{11}(t).
\end{align*}
    Так что, интеграл в \eqref{eq11} сходится равномерно по $0\le x<1$.

    Далее, используем следующее представление для ПФ $G_i(t;x)$ из~\cite{Imomov2012}:
\begin{equation*}
    G_i(t;x)={{qx}\over{i\beta^t}}\cdot{\partial\over{\partial s}}\left[{\left({{{\Phi(t;s)}\over q}}\right)^i}\right]_{s=qx}.
\end{equation*}
    Ввиду последнего равенства в условиях теоремы, с учетом формулы \eqref{eq3},
    функцию $\tau(t;x)$ можно преобразовывать к следующему виду:
\begin{equation}              \label{eq12}
    \tau(t;x)={\Phi'(t;x)-P_{11}(t)},
\end{equation}
    здесь и далее знак производной понимается по переменной $x$.

    Теперь вспомним дифференциальный аналог Основной Леммы теории критических МВС ($f'(1-) = 0$) из \cite{Imomov2005}.
    Делаем замену переменной $x=\exp\{-\lambda\varepsilon(t)\}$ в выражении $\Phi'(t;x)$ так, чтобы $\varepsilon(t)\downarrow 0$,
    $\lambda\to\infty$ и $\lambda=o\left({1}\big/{\varepsilon(t)}\right)$. Тогда очевидно, $x\uparrow 1$ и
\begin{equation}                  \label{eq13}
    \left.{{{\partial\Phi(t;x)}\over{\partial x}}}\right|_{x =\exp\{-\lambda\varepsilon (t)\}}
    =\left[{1-{{1-\Phi\left({t;\exp\{-\lambda\varepsilon (t)\}} \right)}\over{1-\Phi\left({t;0}\right)}}}\right]^2
    \left({1+o(1)}\right)\quad \text{при}\quad t \to \infty.
\end{equation}
    В свою очередь, из многочисленных источников известно, что сама Основная Лемма утверждает справедливости разложения
\begin{equation}           \label{eq14}
    1-\Phi(t;x)={{1-x}\over{1+{{\displaystyle bt}\over\displaystyle 2}\left({1-x}\right)}}
    \left({1+\mathcal{O}\left({{{\ln t}\over t}}\right)}\right)
    \quad \text{при}\quad t \to \infty
\end{equation}
    равномерного для всех $0\le x\leq{r}<1$; см.~{\cite[p.~74]{Sev71}}. Из вида \eqref{eq14} непосредственно убедимся, что
\begin{equation}\label{eq15}
    1-{{1-\Phi(t;x)}\over {1-\Phi(t;0)}}=\mathcal{O}\left({{{\ln t}\over t}}\right)\quad \text{при}\quad t \to \infty
\end{equation}
    для всех $0\le x<1$. Тогда ввиду \eqref{eq13}, при $x\uparrow 1$ следующая оценка очевидна:
\begin{equation}\label{eq16}
    {{\partial\Phi(t;x)}\over{\partial x}}=\mathcal{O}\left({{{\ln^2t}\over{t^2}}}\right)\quad \text{при}\quad t \to \infty.
\end{equation}

    Преобразование $x=\exp\{-\lambda\varepsilon(t)\}$ в левой части (\ref{eq15}) приводит нас к соотношению
\begin{equation*}
    {{\lambda\varepsilon(t)}\over{\lambda\varepsilon(t)+
    {\displaystyle 2\over{\displaystyle bt}}}}=1+o\left(1\right) \quad \text{при}\quad t \to \infty,
\end{equation*}
    чего рассмотрим как асимптотическое уравнение относительно $\varepsilon (t)$, не забывая при этом, что
    $\lambda\to\infty$. Продолжая рассуждение убедимся, что для того чтобы левая часть стремилась к 1, достаточно
    чтобы порядок убывания бесконечно малой величины $\varepsilon (t)$ был $\mathcal{O}\left({1}\big/{t}\right)$.
    Так что, для сохранения закономерности сказанных рассуждений, необходимо положить $\varepsilon (t)={C}\big/{t}$,
    где $C-$любая постоянная. Не нарушая общности, и с желанием упрощения формул, с этого место мы выбираем
    $\varepsilon(t)={2}\big/{bt}$. Тогда, стандартные вычисления нам дает следующее соотношение:
\begin{equation}          \label{eq17}
    \left.{{{\partial\Phi(t;x)}\over{\partial x}}}\right|_{x=\exp\{-\lambda\varepsilon (t)\}}
    ={1\over{\bigl({1+\lambda}\bigr)^2}}\left({1+\mathcal{O}\left({\lambda{{\ln t}\over t}}\right)}\right)
    \quad \text{при}\quad t \to \infty.
\end{equation}

    Выбираем теперь величину $\lambda$. Сравнительный анализ соотношений \eqref{eq16} и \eqref{eq17}
    показывает, что нам необходимо положить  $\lambda=\mathcal{O}\bigl({t}\big/{\ln{t}}\bigr)$.
    Пусть, далее $\lambda={t}\big/{\ln{t}}$.

    Теперь с помощью соотношений \eqref{eq12} и \eqref{eq17}, оценим интеграл в правой части \eqref{eq11}.
    Следуя методу из работы \cite{Nagaev} (см.~также~\cite{Badal}), положим
\begin{equation}            \label{eq18}
    \int_0^1{{{\tau(t;x)}\over{x}}dx}=\left[{\int_{\exp\{-\lambda\varepsilon (t)\}}^1{}+
    \int_0^{\exp\{-\lambda\varepsilon (t)\}}{}}\right]{{\tau(t;x)}\over{x}}dx=:I_1(t)+I_2(t).
\end{equation}
    Делая замену $x=\exp\{-u\varepsilon (t)\}$ в интеграле $I_1(t)$, при $\varepsilon(t)={2}\big/{bt}$, имеем
\begin{align*}
    I_1(t)
    & = \varepsilon (t)\int_0^\lambda{\Bigl[{\Phi'\left({t;\exp\bigl\{-u\varepsilon (t)\bigr\}}\right)-P_{11}(t)}\Bigr]du}    \nonumber\\
\nonumber\\
    & = {2\over bt}\int_0^\lambda{\left[{{1\over{\left({1+u}\right)^2}}
    \left({1+\mathcal{O}\left({u{{\ln t}\over t}}\right)}\right)-P_{11}(t)}\right]du}
    \quad \text{при}\quad t \to \infty     \nonumber.
\end{align*}
    При оценке последнего интеграла учитываем утверждения \eqref{eq2}, \eqref{eq4} и того факта,
    что $\lambda={t}\big/{\ln{t}}$ и, после стандартных аналитических рассуждений, находим
\begin{equation}                \label{eq19}
    I_1(t)={2\over {bt}}\left({1+\mathcal{O}\left({{{\ln^2t}\over t}}\right)}\right)  \quad \text{при}\quad t \to \infty.
\end{equation}

    Чтобы оценить $I_2(t)$, достаточно воспользоваться монотонностью функции $h(t;x)$. Пользуясь соотношением \eqref{eq2},
    \eqref{eq4} и \eqref{eq13}, учитывая при этом соотношений $\varepsilon(t)={2}\big/{bt}$ и $\lambda={t}\big/{\ln{t}}$,
    получаем следующую оценку:
\begin{equation*}
    I_2(t)\le{{{\Phi'}\left({t;\exp\left\{-{2}\big/{b\ln{t}}\right\}}\right)-P_{11}(t)}
    \over {\exp\bigl\{-{2}\big/{b\ln{t}}\bigr\}}}
    =\mathcal{O}\left({{{\ln^2t}\over {t^2}}}\right)    \quad \text{при}\quad t \to \infty.
\end{equation*}
    Рассмотрев последнюю оценку вместе с равенствами \eqref{eq11}, \eqref{eq18} и оценкой \eqref{eq19}, получаем
\begin{equation}\label{eq20}
    \Sigma_1 (t)={2\over t} \left({1+\mathcal{O}\left({{{\ln^2t}\over t}}\right)}\right)    \quad \text{при}\quad t \to \infty.
\end{equation}

    Приступим к оценке второй суммы $\Sigma_2 (t)$ в \eqref{eq10}. Используя рассуждение из равенств \eqref{eq11}, находим
\begin{equation}          \label{eq21}
    \Sigma_2(t)=b\sum_{k\in{\mathcal{E}}}{{{\mathcal{Q}_{1k+1}(t)}\over {k^2}}} = b\int_0^1{{{\mu(t;x)}\over x}dx,}
\end{equation}
    здесь
\begin{equation*}
    \mu(t;x)= {\sum_{k\in{\mathcal{E}}}{{{\mathcal{Q}_{1k+1}(t)}\over k}}x^k}.
\end{equation*}
    Используя опять рассуждение из равенств \eqref{eq11} и, учитывая теперь тот факт,
    что $\int_{0}^{x}s^{k-1}ds={{x^{k}}\big/{k}}$, последнее равенство преобразуем к виду
\begin{equation*}
    \mu(t;x)= \int_0^x{{{\tau(t;s)}\over{s}}ds},
\end{equation*}
    где как и прежде,
\begin{equation*}
    \tau(t;x)= {{G(t;x)-x\mathcal{Q}_{11}(t)}\over{x}}\raise0.9pt\hbox{.}
\end{equation*}

    Подынтегральная функция ${\mu(t;x)}\big/{x}$ в правой части равенства \eqref{eq21} 
    имеет конечные значения в концах области своего определения. Действительно, используя правило Лопиталя вычислим, что
\begin{align*}
    \lim_{x\downarrow 0}{{\mu(t;x)}\over x} &=  \lim_{x\downarrow 0} {{\tau(t;x)}\over{x}}=\mathcal{Q}_{12}(t), \\
\\  \nonumber
    \lim_{x\uparrow 1}{{\mu(t;x)}\over x} &= \lim_{x\uparrow 1}\int_0^x {{{\tau(t;x)}\over{x}}ds}=\mathcal{O}\Bigl({\Sigma_1(t)}\Bigr).
\end{align*}
    Кроме этого, в силу монотонности функции ${\tau(t;x)}\big/{x}$
\begin{equation*}
    \tau(t;x)-\mu(t;x)=\tau(t;x)-\int_0^x {{{\tau(t;s)}\over{s}}ds}\ge 0.
\end{equation*}
    Следовательно
\begin{equation*}
    {\partial\over{\partial x}}\left[{{{\mu(t;x)}\over x}}\right]={{\tau(t;x)-\mu(t;x)}\over{x^2}}\ge 0.
\end{equation*}
    Последнее утверждает, что функция ${\mu(t;x)}\big/ {x}$ монотонно возрастает. Она же ограничена
    в области интегрирования. Отсюда следует что, интеграл в правой части \eqref{eq21} сходится.
    Оценим этот интеграл. Ссылаясь на формулу \eqref{eq3}, легко находим
\begin{align*}
    \mu(t;x)
    & = \sum_{k\in{\mathcal{E}}}P_{1k+1}(t)x^{k}+
    \sum_{k\in{\mathcal{E}}}{{{P_{1k+1}(t)}\over {k}}}x^{k}    \nonumber\\
\nonumber\\
    & = \frac{\Phi(t;x)-P_{10}(t)-P_{11}(t)x}{x}+
    \int_{0}^{x}\frac{\Phi(t;s)-P_{10}(t)-P_{11}(t)s}{s^{2}}ds.   \nonumber
\end{align*}
    Следовательно,
\begin{equation}          \label{eq22}
    \Sigma_2(t)=b\int_{0}^{1}\frac{v(t;x)}{x}dx +
    b\int_{0}^{1}\frac{1}{x}\left[\int_{0}^{x}\frac{v(t;s)}{s}ds\right]dx=:J_{1}(t)+J_{2}(t),
\end{equation}
    где
\begin{equation*}
    v(t;x)=\frac{\Phi(t;x)-P_{10}(t)-P_{11}(t)x}{x}.
\end{equation*}

    Итак, оценке подлежат интегралы $J_{1}(t)$ и $J_{2}(t)$. При этом мы по-прежнему придерживаемся замены
    переменной $x=\exp\bigl\{-u\varepsilon(t)\bigl\}$. Для нашей ближайшей цели нам понадобится следующее асимптотическое
    соотношение, которого можно получить с помощью разложения (\ref{eq14}), сохраняя при этом прежние обозначения:
\begin{equation}         \label{eq23}
    v\left(t;e^{-\lambda\varepsilon(t)}\right)=\frac{1}{1+\lambda}{2\over bt}
    \left({1+\mathcal{O}\left({\lambda{{\ln t}\over t}}\right)}\right)-P_{11}(t)   \quad \text{при}\quad t \to \infty.
\end{equation}
    Элементарное рассуждение показывает, что функция $v(t;x)\big/x$ монотонно возрастает
    по $0\leq{x}<1$, причем она имеет конечные значения в концах этой области:
\begin{align*}
    \lim_{x\downarrow 0}{{v(t;x)}\over x} &=  P_{12}(t), \\
\\  \nonumber
    \lim_{x\uparrow 1}{{v(t;x)}\over x}  &=  1-P_{10}(t)-P_{11}(t).
\end{align*}
    Поэтому интеграл $J_{1}(t)$ в правой части \eqref{eq22} сходится
    равномерно по $0\leq x<1$. Следуя предыдущему рассуждению, положим
\begin{equation}      \label{eq24}
    J_{1}(t)=b\left[{\int_{\exp\{-\lambda\varepsilon (t)\}}^1{}+\int_0^{\exp\{-\lambda\varepsilon (t)\}}}\right]
    {{v(t;x)}\over x}dx=:J_{11}(t)+J_{12}(t).
\end{equation}

    Делая, опять, замену $x=\exp\{-u\varepsilon (t)\}$ теперь в интеграле $J_{11}(t)$, при этом принимая во
    внимание \eqref{eq23} и $\varepsilon(t)={2}\big/{bt}$, получаем следующее асимптотическое соотношение:
\begin{align*}
    J_{11}(t)
    & = {2\over t} \int_0^\lambda {v\left({t;e^{-u\varepsilon (t)}}\right)du}    \nonumber\\
\nonumber\\
    & = {2\over t} \int_0^\lambda {\left[\frac{1}{1+u} {2\over bt}\left({1+\mathcal{O}\left({u{{\ln t}
    \over t}}\right)}\right)-P_{11}(t)\right]du}  \quad \text{при}\quad t \to \infty.   \nonumber
\end{align*}
    Отсюда, используя формулу \eqref{eq3} и учитывая тот факт, что
    $\lambda={t}\big/{\ln{t}}$, с учетом  утверждения \eqref{eq4}, находим
\begin{equation}                \label{eq25}
    J_{11}(t)=\mathcal{O}\left({{{\ln t}\over t^{2}}}\right)  \quad \text{при}\quad t \to \infty.
\end{equation}
    В свою очередь, ввиду монотонности функции $v(t;x)\big/x$ и соотношения \eqref{eq23},
\begin{align}                \label{eq26}
    J_{12}(t)
    & = b{\int_0^{\exp\{-\lambda\varepsilon (t)\}}}{{v(t;x)}\over x}dx    \nonumber\\
\nonumber\\
    & \leq bv\left(t;e^{-\lambda\varepsilon(t)}\right)=
    \mathcal{O}\left({{{\ln t}\over t^{2}}}\right)    \quad \text{при}\quad t \to \infty.
\end{align}
    Из \eqref{eq24}--\eqref{eq26}, получим оценку
\begin{equation}                 \label{eq27}
    J_{1}(t)=\mathcal{O}\left({{{\ln t}\over t^{2}}}\right)     \quad \text{при}\quad t \to \infty.
\end{equation}

    Займемся теперь оценкой интеграла $J_{2}(t)$. Запишем его в виде
\begin{equation*}
    J_{2}(t)=b\int_{0}^{1}\frac{1}{x}\left[\int_{0}^{x}\frac{v(t;s)}{s}ds\right]dx
    =:b\int_{0}^{1}\frac{\omega(t;x)}{x}dx.
\end{equation*}
    Используя правило Лопиталя, находим
\begin{align*}
    \lim_{x\downarrow 0}{{\omega(t;x)}\over x} &=  \lim_{x\downarrow 0}{{v(t;x)}\over x}=P_{12}(t), \\
\\  \nonumber
    \lim_{x\uparrow 1}{{\omega(t;x)}\over x}  &=  \lim_{x\uparrow 1}\int_{0}^{x}{{v(t;s)}\over s}ds=J_{1}(t).
\end{align*}
    Так что исследуемый интеграл сходится по $x$. Пусть
\begin{equation}            \label{eq28}
    J_{2}(t)=b\left[{\int_{\exp\{-\lambda\varepsilon (t)\}}^1{} +
    \int_0^{\exp\{-\lambda\varepsilon (t)\}}}\right] {{\omega(t;x)}\over x}dx=:J_{21}(t)+J_{22}(t).
\end{equation}

    Заменой $x=\exp\bigl\{-u\varepsilon(t)\bigl\}$ в интеграле $J_{21}(t)$, с учетом
    монотонности функции $v(t;x)\big/x$, мы получим следующую цепочку соотношений:
\begin{align}                \label{eq29}
    J_{21}(t)
    & = b{\int_{\exp\{-\lambda\varepsilon (t)\}}^{1}}{{\omega(t;x)}\over x}dx
    =\frac{2}{t} \int_{0}^{\lambda}{\omega\left(t;e^{-u\varepsilon (t)}\right)}du    \nonumber\\
\nonumber\\
    & = \frac{2}{t} \int_{0}^{\lambda}\left[\int_{0}^{\exp\{-u\varepsilon (t)\}}{{v(t;x)}\over x}dx\right]du
    \leq \frac{2}{t} \int_{0}^{\lambda}{v\left(t;e^{-u\varepsilon (t)}\right)}du=J_{11}(t).
\end{align}

    Далее, свойства же монотонности функции $v(t;x)\big/x$ нам дает
\begin{equation*}
    v(t;x)-\omega(t;x)=v(t;x)-\int_{0}^{x}{{v(t;s)}\over s}ds\geq 0.
\end{equation*}
    Следовательно,
\begin{equation*}
    \frac{\partial}{\partial x}\left[{{\omega(t;x)}\over x}\right]=\frac{v(t;x)-\omega(t;x)}{x^{2}}\geq 0,
\end{equation*}
    что свидетельствует о монотонности функции $\omega(t;x)\big/x$. Используя это свойство,
    нетрудно оценить второй интеграл в равенстве \eqref{eq28}:
\begin{align}                \label{eq30}
    J_{22}(t)
    & = b{\int_{0}^{\exp\{-\lambda\varepsilon (t)\}}}{{\omega(t;x)}\over x}dx    \nonumber\\
\nonumber\\
    & \leq b\omega\left(t;e^{-u\varepsilon(t)}\right)=b{\int_{0}^{\exp\{-\lambda\varepsilon (t)\}}}{{v(t;x)}\over x}dx=J_{12}(t).
\end{align}
    Рассмотрев вместе \eqref{eq28}--\eqref{eq30}, с учетом \eqref{eq24}, заключаем
\begin{equation}             \label{eq31}
    J_{2}(t)=\mathcal{O}\bigl(J_{1}(t)\bigr)   \quad \text{при}\quad t \to \infty.
\end{equation}

    Наконец, собирая вместе равенства \eqref{eq10}, \eqref{eq22} и оценки \eqref{eq20}, \eqref{eq27}, \eqref{eq31},
    мы получаем требуемое утверждение, что завершает доказательство Теоремы~\ref{INTh:1}.
\end{proof}

\begin{proof}[\textbf{Доказательство Теоремы~\ref{INTh:2}}]
    Рассмотрим преобразования Харриса-Севастьянова
\begin{equation*}
    f_{q}(x):={{f\left({qx}\right)}\over q} \qquad \mbox{и}
    \qquad \Phi_{q}(t;x):={{\Phi\left({t;qx}\right)}\over q}\raise0.9pt\hbox{.}
\end{equation*}
    Нетрудно заметить, что $f_{q}(x)$ представляет собой инфинитезимальную ПФ, которая, в свою очередь, порождает докритическую
    МВС $\left\{Z_q(t)\right\}$ с средним значением интенсивностей ${f'_{q}(1-)}=\ln{\beta}$ и пространством возможных состояний
    $\mathcal{S}_{0}=\left\{0\right\}\cup{\mathcal{S}}$. А ПФ распределения потомков одной частицы этой системы
\begin{equation*}
    {\textsf{E}}\left[{x^{Z_q(t)}\Bigm|{Z_q(0)=1}}\right]
    =\sum_{j\in{\mathcal{S}_{0}}}{P_{1j}^{(q)}(t)x^j} = \Phi_{q}(t;x),
\end{equation*}
    где $P_{1j}^{(q)}(t)=q^{j-1} P_{1j}(t)$ и $P_{1j}(t)$ -- распределение вероятностей
    числа поколений одной частицы в МВС $\left\{Z(t)\right\}$, порожденной инфинитезимальной ПФ $f(x)$.

    В работе \cite[Теорема~7]{Imomov2012} доказано, что если $f'(1-)>0$, то
\begin{equation*}
    G_i (t;x) \longrightarrow x\exp\left\{{\int_{qx}^q {{{f'(s)-f'(q)}\over{f(s)}}ds}}\right\}
    =: U(x) \quad \text{при}\quad t \to \infty
\end{equation*}
    для всех $i\in{\mathcal{E}}$, где предельная ПФ $U(x)$ порождает инвариантную меру для надкритических МВС.
    Ссылаясь на \cite{ImomovTSU2021}, мы убедимся в том, что сходимость $G_i(t;x) \to U(x)$ верна и
    для системы $\left\{Z_q(t)\right\}$. Следовательно, в наших обозначениях, при $i=1$ справедливо
    следующее асимптотическое соотношение:
\begin{equation}         \label{eq32}
    G(t;x) = x\kappa(x)\bigl({1+o(1)}\bigr)    \quad \text{при}\quad t \to \infty,
\end{equation}
    где
\begin{equation*}
    \kappa(x)= \exp\left\{{\int_x^1{{{f'_{q}(s)-f'_{q}(1-)}\over{f_{q}(s)}}ds}}\right\}.
\end{equation*}
    В свою очередь, используя формулу Тейлора вблизи точки $x=1$, с учетом тот факт, что $f''_{q}(1-)=b$,
    мы легко находим следующее локально-асимптотическое представление:
\begin{equation*}
    \kappa(x)\sim 1-\gamma\left({1-x}\right) \quad \text{при}\quad  x\uparrow 1,
\end{equation*}
    здесь по-прежнему $\gamma:={b\big/ \bigl|{\ln\beta}\bigr|}$. Отсюда, полагая
    $x=\exp\bigl\{-\lambda\beta^t\bigr\}$ так, чтобы $\lambda\to\infty$ и $\lambda\beta^t =o(1)$,
    получаем теперь следующую асимптотическую формулу:
\begin{equation}              \label{eq33}
    \kappa\left({e^{-\lambda\beta^t}}\right)\sim 1-\gamma\lambda\beta^t \quad \text{при}\quad  t\to\infty.
\end{equation}

    Теперь, из \eqref{eq8} и \eqref{eq9} для случая $\beta<1$, получим
\begin{align}                 \label{eq34}
    {\textsf{Var}}\widehat\beta(t)
    & = \left({1+\gamma}\right)\beta\left({1-\beta}\right)\sum_{k\in{\mathcal{E}}}{{{\mathcal{Q}_{1k+1}(t)} \over k}}  \nonumber\\
\nonumber\\
    & \quad +\gamma\left({1-\beta}\right)\sum_{k\in{\mathcal{E}}}{{{\mathcal{Q}_{1k+1}(t)}\over{k^2}}}
    =:\Sigma_{\gamma{1}} (t)+\Sigma_{\gamma{2}} (t).
\end{align}
    Далее мы следуем методу доказательства Теоремы~\ref{INTh:1}. Первую сумму в \eqref{eq34} запишем виде
\begin{equation}             \label{eq35}
    \Sigma_{\gamma{1}}(t)=\left({1+\gamma}\right)\beta\left({1-\beta}\right)\int_0^1{{{\tau(t;x)}\over{x}}dx},
\end{equation}
    где по-прежнему
\begin{equation*}
    \tau(t;x)={{G(t;x)-x\mathcal{Q}_{11}(t)}\over{x}}\raise0.9pt\hbox{.}
\end{equation*}
    Имеем
\begin{equation}             \label{eq36}
    \int_0^1 {{{\tau(t;x)}\over{x}}dx}=\left[{\int_{\exp\left\{-\lambda\beta^t\right\}}^{1} +
    \int_0^{\exp\left\{-\lambda\beta^t\right\}}}\right]{{\tau(t;x)}\over{x}}dx
    =:\mathcal{I}_{\tau{1}}(t)+\mathcal{I}_{\tau{2}}(t),
\end{equation}
    где $\lambda=o\bigl({\beta^{-t}}\bigr)\to\infty$.

    Делая замену $x =\exp\left\{-s\beta^t\right\}$ в интеграле $\mathcal{I}_{\tau{1}}(t)$ и,
    используя соотношений \eqref{eq32} и \eqref{eq33}, получим
\begin{align*}
    \mathcal{I}_{\tau{1}}(t)
    & = \int_{\exp\{-\lambda\beta^t\}}^1{{\tau(t;x)}\over{x}}dx
        =\beta^t\int_0^\lambda{\left[\kappa\left({e^{-\lambda\beta^t}}\right)-\mathcal{Q}_{11}(t)\right]}ds    \nonumber\\
\nonumber\\
    & = \beta^t\int_0^\lambda{\Bigl[{1-s\gamma\beta^t\bigl({1+o(1)}\bigr)-\mathcal{Q}_{11}(t)}\Bigr]}ds \quad \text{при}\quad  t\to\infty.
\end{align*}
    Отсюда, в силу утверждения \eqref{eq5} и $\lambda=o\bigl({\beta^{-t}}\bigr)$, получаем оценку
\begin{equation}             \label{eq37}
    \mathcal{I}_{\tau{1}}(t)=o(1) \quad \text{при}\quad  t\to\infty.
\end{equation}

    Чтобы оценить $\mathcal{I}_{\tau{2}}(t)$ вспомним, что функция ${\tau(t;x)}\big/{x}$ монотонно
    возрастает. Тогда, учитывая опять утверждение \eqref{eq5} и $\lambda=o\bigl({\beta^{-t}}\bigr)$, имеем
\begin{align}                \label{eq38}
    0<\mathcal{I}_{\tau{2}}(t)
    & \leq  \tau\left({t;e^{-\lambda\beta^t}}\right)   \nonumber\\
\nonumber\\
    & = 1-\mathcal{Q}_{11}(t)-\lambda\gamma\beta^t \bigl({1+o(1)}\bigr)
    \quad \text{при}\quad  t\to\infty.
\end{align}
    Собирая равенства \eqref{eq35} и \eqref{eq36} вместе с оценками \eqref{eq37} и \eqref{eq38}, находим
\begin{equation*}
    \Sigma_{\gamma{1}}(t)=\left({1+\gamma}\right)\beta\left({1-\beta}\right)
    \bigl(1-\mathcal{Q}_{11}(t)\bigr)+o(1)  \quad \text{при}\quad  t\to\infty.
\end{equation*}
    Сходимость \eqref{eq5} утверждает, что $\mathcal{Q}_{11}(t)$ имеет зависящий от $q$ конечный предел. Поэтому
\begin{equation}           \label{eq39}
    \Sigma_{\gamma{1}}(t)=\mathcal{O}\left(1\right) \quad \text{при}\quad  t\to\infty.
\end{equation}

    Переходим к оценке суммы $\Sigma_{\gamma{2}}(t)$. Для этого следуем за ходом рассуждения соответствующей
    части доказательства Теоремы~\ref{INTh:1}, начиная с формулы \eqref{eq21}. Имеем
\begin{align}                \label{eq40}
    \Sigma_{\gamma{2}}(t)
    & =  \gamma\left({1-\beta}\right) \left[{\int_{\exp\left\{-\lambda\beta^t\right\}}^{1} 
        + \int_0^{\exp\left\{-\lambda\beta^t\right\}}}\right]{{\mu(t;x)}\over{x}}dx   \nonumber\\
\nonumber\\
    & =: \gamma\left({1-\beta}\right) \Bigl[\mathcal{J}_{\mu{1}}(t)+\mathcal{J}_{\mu{2}}(t)\Bigr],
\end{align}
    где $\lambda=o\bigl({\beta^{-t}}\bigr)\to\infty$ при $t\to\infty$.

    Оценим интеграл $\mathcal{J}_{\mu{1}}(t)$. Делаем в нем замену $x =\exp\left\{-u\beta^t\right\}$ и,
    используя монотонности функции ${\tau(t;x)}\big/{x}$ и соотношений \eqref{eq32} и \eqref{eq33},
    получим следующую цепочку соотношений:
\begin{align}                \label{eq41}
    0<\mathcal{I}_{\mu{1}}(t)
    & = \beta^{t} \int_{0}^{\lambda}{\mu\left(t;e^{-u\beta^{t}}\right)}du    \nonumber\\
\nonumber\\
    & = \beta^{t} \int_{0}^{\lambda}
    \left[\int_{0}^{\exp\left\{-u\beta^{t}\right\}}{{{\tau(t;s)}\over s}ds}\right]du   \nonumber\\
\nonumber\\
    & \leq \beta^{t} \int_{0}^{\lambda}{\tau\left(t;e^{-u\beta^{t}}\right)}du
    =\beta^{t} \int_0^\lambda{\left[\kappa\left({e^{-\lambda\beta^t}}\right)-\mathcal{Q}_{11}(t)\right]}ds  \nonumber\\
\nonumber\\
    & =\Bigl[1-\mathcal{Q}_{11}(t)-\gamma\beta^t\bigl({1+o(1)}\bigr)\Bigr]\lambda\beta^{t}  \quad \text{при}\quad  t\to\infty.
\end{align}
    Поскольку $\lambda\beta^t =o(1)$, то в силу \eqref{eq5} из \eqref{eq41} получаем оценку
\begin{equation}           \label{eq42}
    \mathcal{I}_{\mu{1}}(t)={o}\left(1\right) \quad \text{при}\quad  t\to\infty.
\end{equation}

    Для оценки интеграла $\mathcal{I}_{\mu{2}}(t)$ последовательно используем свойство монотонности
    функций ${\mu(t;x)}\big/{x}$ и ${\tau(t;x)}\big/{x}$. В результате, согласно \eqref{eq5}
    получаем следующую цепочку соотношений:
\begin{align*}
    0 < \mathcal{I}_{\mu{2}}(t)
    & \leq \mu\left({t;e^{-\lambda\beta^t}}\right)=\int_0^{\exp\left\{-\lambda\beta^t\right\}}{{\tau(t;s)}\over{s}}ds \nonumber\\
\nonumber\\
    & \leq \tau\left({t;e^{-\lambda\beta^t}}\right)=\kappa\left({e^{-\lambda\beta^t}}\right)-\mathcal{Q}_{11}(t)   \nonumber\\
\nonumber\\
    & = 1-\mathcal{Q}_{11}(t)-\gamma\lambda\beta^t \bigl({1+o(1)}\bigr) \longrightarrow 1-{{\left|{\ln\beta}\right|\mathcal{A}}\over{a_0}}
    \quad \text{при}\quad  t\to\infty.
\end{align*}
    Таким образом
\begin{equation}           \label{eq43}
    \mathcal{I}_{\mu{2}}(t) =\mathcal{O}\left(1\right) \quad \text{при}\quad  t\to\infty.
\end{equation}

    Требуемое утверждение вытекает теперь из комбинации соотношений \eqref{eq34}, \eqref{eq39}, \eqref{eq40}, \eqref{eq42} и \eqref{eq43}.

    Теорема~\ref{INTh:2} доказана.
\end{proof}

\label{lastpage}

\end{document}